\documentclass[final,1p,times]{elsarticle}

\usepackage{amssymb}
 \usepackage{amsthm}
\usepackage{amscd}
\usepackage{amsmath}
\usepackage{amsfonts}
\usepackage{amssymb}
\usepackage{graphicx}
\newtheorem{theorem}{Theorem}

\newtheorem{lemma}[theorem]{Lemma}

\usepackage{mathrsfs}
\usepackage{titletoc}

\newcommand{\ra}{\rightarrow}

\newcommand{\f}{\frac}

\newcommand{\be}{\begin{equation}}
\renewcommand{\ra}{\rightarrow}
\newcommand{\ee}{\end{equation}}
\newcommand{\bea}{\begin{eqnarray}}
\newcommand{\eea}{\end{eqnarray}}
\newcommand{\bna}{\begin{eqnarray*}}
\newcommand{\ena}{\end{eqnarray*}}

\renewcommand{\le}{\left}
\newcommand{\ri}{\right}

\usepackage{color}

\journal{***}

\begin{document}

\begin{frontmatter}

\title{Extremal functions for Adams' inequalities in dimension four}

\author{Xiaomeng Li$^{1,2}$}
\ead{xmlimath@ruc.edu.cn}
\address{$^1$ School of Information, Huaibei Normal University, Huaibei, 235000, P. R. China}
\address{$^2$Department of Mathematics,
Renmin University of China, Beijing 100872, P. R. China}

\begin{abstract}
Let $\Omega\subset \mathbb{R}^4$ be a smooth bounded domain, $W_0^{2,2}(\Omega)$ be the usual Sobolev space. For any positive integer $\ell$, $\lambda_{\ell}(\Omega)$ is the $\ell$-th eigenvalue of the bi-Laplacian operator. Define $E_{\ell}=E_{\lambda_1(\Omega)}\oplus E_{\lambda_2(\Omega)}\oplus\cdots\oplus E_{\lambda_{\ell}(\Omega)}$, where $E_{\lambda_i(\Omega)}$ is eigenfunction space associated with $\lambda_i(\Omega)$. $E^{\bot}_{\ell}$ denotes the orthogonal complement of $E_\ell$ in $W_0^{2,2}(\Omega)$. For $0\leq\alpha<\lambda_{\ell+1}(\Omega)$, we define a norm by $\|u\|_{2,\alpha}^{2}=\|\Delta u\|^2_2-\alpha \|u\|^2_2$ for $u\in E^\bot_{\ell}$. In this paper, using the blow-up analysis, we prove the following Adams inequalities
$$\sup_{u\in E_{\ell}^{\bot},\,\| u\|_{2,\alpha}\leq 1}\int_{\Omega}e^{32\pi^2u^2}dx<+\infty;$$
moreover, the above supremum can be attained by a function $u_0\in E_{\ell}^{\bot}\cap C^4(\overline{\Omega})$ with $\|u_0\|_{2,\alpha}=1$.
 This result extends that of Yang (J. Differential Equations, 2015), and complements that of Lu and Yang (Adv. Math. 2009) and Nguyen (arXiv: 1701.08249, 2017).

\end{abstract}


\end{frontmatter}
\section{Introduction}

Trudinger-Moser inequalities play important roles in  analysis and geometry.
There are two interesting subjects in the study of Trudinger-Moser inequalities: one is the attainability of the best constants, the other is the existence of extremal functions. The research on sharp constants was initiated by Yudovich \cite{Yudovich}, Pohozaev \cite{19} and Trudinger \cite{22}. Later Moser \cite{14} found the best constant: if $\beta\leq\beta_0=n\omega_{n-1}^{1/(n-1)}$, then
\be\label{Moser}\sup_{u\in W_0^{1,n}(\Omega),\,\|\nabla u\|_n=1}\int_{\Omega}e^{\beta|u|^{n/(n-1)}}dx<\infty,\ee
where $\Omega$ is an open subset of $\mathbb{R}^n$ $(n\geq2)$ with finite Lebesgue measure, $\omega_{n-1}$ is the measure of the unit sphere in $\mathbb{R}^n$;
moreover, if $\beta>\beta_0$, the integrals in (\ref{Moser}) are still finite, but the supremum is infinite.
The sharp constants for higher order derivatives of Moser's inequality was due to Adams \cite{Adams}. For any fixed positive integer $m$, let $u\in C_0^m(\Omega)$, the space of functions having $m$-th continuous derivatives and compact support. To state Adams' result, we use the symbol $\nabla^m u$ to denote the $m$-th order gradient for $u$. Precisely
\bna\nabla^m u=\le\{\begin{array}{lll}\Delta^{\frac{m}{2}}u& {\rm when} \ m \ {\rm is} \ {\rm even},\\[1.2ex]
\nabla\Delta^{\frac{m-1}{2}}u& {\rm when} \ m \ {\rm is} \ {\rm odd},
\end{array}\ri.\ena
where $\nabla$ and $\Delta$ denote the usual gradient and the Laplacian operators.  Adams proved that if $\beta\leq\beta(n,m)$ and $0<m<n$, then
\be\label{AD}
\sup_{u\in W_0^{m,\frac{n}{m}}(\Omega),\,\|\nabla^m u\|_{L^{\frac{n}{m}}(\Omega)}\leq1}\int_{\Omega}e^{\beta |u|^{n/(n-m)}}dx\leq C_{m,n}|\Omega|
\ee
for some constant $C_{m,n}$,
where
\bna\beta(n,m)=\le\{\begin{array}{lll}\frac{n}{\omega_{n-1}}\le[\frac{\pi^{n/2}2^m\Gamma(\frac{m+1}{2})}{\Gamma(\frac{n-m+1}{2})}\ri]^{\frac{n}{n-m}} & {\rm when} \ m \ {\rm is} \ {\rm odd},\\[1.2ex]
\frac{n}{\omega_{n-1}}\le[\frac{\pi^{n/2}2^m\Gamma(\frac{m}{2})}{\Gamma(\frac{n-m}{2})}\ri]^{\frac{n}{n-m}} & {\rm when} \ m \ {\rm is} \ {\rm even}.
\end{array}\ri.\ena
Moreover, $\beta(n,m)$ is the best constant in the sense that if $\beta>\beta(n,m)$, then the supremum in (\ref{AD}) is infinite.
The manifold version of Adams' inequality was obtained by Fontana \cite{Fontana}.

Extremal functions for (\ref{Moser}) were first found by Carleson and Chang \cite{C-C} when $\Omega$ is the unit ball in $\mathbb{R}^n$.
This result was then extended by Flucher \cite{Flucher} to a general domain $\Omega\subset \mathbb{R}^2$, and by Lin \cite{Lin} to a bounded smooth domain $\Omega\subset\mathbb{R}^n$ $(n\geq 2)$.
  In $2004$, it was proved by Adimurthi and Druet \cite{A-O} that for any $\alpha$, $0\leq\alpha<\lambda_1(\Omega)$, there holds
\be\label{ao}\sup_{u\in W_0^{1,2}(\Omega),\,\|\nabla u\|_2\leq1}\int_{\Omega}e^{4\pi^2(1+\alpha\|u\|^2_2)}dx<+\infty\ee
and the supremum is infinit for $\alpha\geq \lambda_1(\Omega)$, where $\lambda_1(\Omega)$ is the first eigenvalue of the Laplacian operator with respect to Dirichlet boundary condition. The inequality (\ref{ao}) was generalized by Yang \cite{Yang2006} to high dimension, by Lu and Yang \cite{Lu-Yang-p} and J. Zhu \cite{Zhu} to the versions involving  $L^p$ norms, by Souza and do \'{O} \cite{DD, DD1} and Ruf \cite{Ruf} and Li and Ruf \cite{Li-Ruf} to the whole Euclidean space, by Tintarev \cite{Tintarev} and Yang \cite{YangJDE-2015} to the following form:
\be\label{Sharpform}\sup_{u\in W_0^{1,2}(\Omega),\,\|\nabla u\|^2_2-\alpha \|u\|^2_2\leq1}\int_{\Omega}e^{4\pi u^2}dx<+\infty\ee
for any $0\leq\alpha<\lambda_{1}(\Omega)$. In particular, Yang \cite{YangJDE-2015} proved that for $0\leq\alpha<\lambda_{\ell+1}(\Omega)$, the extremal function for (\ref{Sharpform}) exists, where $\ell$ is a positive integer and $\lambda_{\ell}$ denotes the $\ell$-th eigenvalue of the Laplacian operator with the Dirichlet boundary condition. The singular version of (\ref{Sharpform}) was considered by Yang and Zhu \cite{Yang-Zhu}, Li and Yang \cite{Li-Yang}, and the author \cite{Li} in $\mathbb{R}^n$ ($n\geq2$). For other works on singular Trudinger-Moser inequalities, we refer the reader to \cite{Adimurthi-Sandeep, Adi-Yang, C-R, M-S} and the references therein.

The study of Turdinger-Moser inequalities on Riemannian manifolds was initiated by Aubin \cite{Aubin} and Cherrier \cite{Cherrier1, Cherrier2}. Much work has also been done in this direction, see for examples \cite{DJLW1997, Li2001, Li2005, Yang-Tran, Yangjfa, YangJGA}.

Let us come back to the Adams inequality in dimension four. Namely
\be\label{Adams}
\sup_{u\in W_0^{2,2}(\Omega),\,\|\Delta u\|_2\leq1}\int_{\Omega}e^{32\pi^2 u^2}dx<+\infty,\ee
where $\Omega\subset \mathbb{R}^4$ is a smooth bounded domain. For any $\alpha$: $0\leq\alpha<\lambda_1(\Omega)$,
it was proved by Lu and Yang \cite{Lu-Yang} that
\be\label{ly}
\sup_{u\in W_0^{2,2}(\Omega),\,\|\Delta u\|^2_2=1}\int_{\Omega}e^{32\pi^2u^2(1+\alpha\|u\|^2_2)}dx<+\infty
\ee
and the supremum is infinite when $\alpha\geq\lambda_1(\Omega)$. Here, by definition,
$$\lambda_1(\Omega)=\inf_{u\in W_0^{2,2}(\Omega),u\not\equiv0}\frac{\|\Delta u\|^2_2}{\|u\|^2_2}.$$
The extremal function for supremum (\ref{ly}) was obtained for $\alpha$ sufficiently small.
This result is recently strengthened by Nguyen \cite{Nguyen} to the following form:
\be\label{N}
\sup_{u\in W_0^{2,2}(\Omega),\,\|\Delta u\|_{2}^2-\alpha \|u\|_2^2\leq1}\int_{\Omega}e^{32\pi^2u^2}dx<+\infty
\ee
for $0\leq\alpha<\lambda_1(\Omega)$, and the above supremum can be achieved by applying the blow-up analysis method.
Motivated by the work \cite{YangJDE-2015}, we shall improve (\ref{N}) to the case involving higher order eigenvalues.
Note that the Dirichlet boundary problem
\bna\le\{\begin{array}{lll}\Delta^2 u=\lambda u \,&{\rm in} \ \Omega,\\[1.2ex]
u=\frac{\partial u}{\partial \nu}=0 \,&{\rm on} \ \partial\Omega
\end{array}\ri.\ena
possesses a sequence of eigenvalues $0<\lambda_1(\Omega)<\cdots<\lambda_i(\Omega)<\lambda_{i+1}(\Omega)<\cdots$. It is known that $\lambda_i(\Omega)\ra \infty$ as $i\ra\infty$, see for example (\cite{Brezis}, Section 6.3).
The corresponding eigenfunction space can be written as
$$E_{\lambda_i(\Omega)}=\le\{u\in W_0^{2,2}(\Omega):\Delta^2 u=\lambda_{i}(\Omega)u\ri\}.$$
For any positive integer $\ell$, we set
$$E_{\ell}=E_{\lambda_1(\Omega)}\oplus E_{\lambda_2(\Omega)}\oplus\cdots\oplus E_{\lambda_{\lambda_{\ell}}(\Omega)}$$
and
\be\label{E}E_{\ell}^\bot=\le\{u\in W_0^{2,2}(\Omega):\int_{\Omega}uvdx=0, \forall v\in E_{\ell}\ri\}.\ee
Clearly  $W_0^{2,2}(\Omega)$ is a real Hilbert space when it is equipped with the inner product
$$<u,v>=\int_{\Omega}\Delta u \Delta vdx, \quad\forall u,v\in W_0^{2,2}(\Omega).$$
According to (\cite{Brezis}, Theorem 9.31), each eigenfunction space $E_{\lambda_i(\Omega)}$ has finite dimension. Suppose
$\dim E_{\lambda_i(\Omega)}=n_i$
and $(e_{ij})$ $(1\leq j\leq n_i,1\leq i\leq \ell)$ be the basis of $E_\ell$. Then
\be\label{ei}\le\{\begin{array}{lll}E_{\lambda_{i} (\Omega)}={\rm span}\{e_{i1},\cdots,e_{in_i}\},\quad i=1,2,\cdots,\ell,\\[1.2ex]
E_\ell={\rm span}\{e_{11},\cdots,e_{1 n_1},e_{21},\cdots,e_{2 n_2},\cdots,e_{\ell1},\cdots,e_{\ell n_\ell}\},\\[1.2ex]
\int_{\Omega}e_{ij}e_{kl}dx=0\quad i\neq k\, {\rm or} \,j\neq l,\\[1.2ex]
\int_{\Omega}e_{ij}^2dx=1.
\end{array}\ri.\ee
Similar as in \cite{YangJDE-2015}, we define
\be\label{lambdal}\lambda_{\ell+1}(\Omega)=\inf_{u\in W_0^{2,2}(\Omega),\,u\in E_{\ell}^\bot,\,u\not\equiv0}\frac{\|\Delta u\|^2_2}{\|u\|^2_2}.\ee
If $u\in E^{\bot}_{\ell}$ satisfies $\int_\Omega|\Delta u|^2dx-\alpha\int_\Omega u^2dx\geq0$, then we denote
\be\label{1alpha}\|u\|_{2,\,\alpha}=\le(\int_\Omega|\Delta u|^2dx-\alpha\int_\Omega u^2dx\ri)^{1/2}.\ee

In this paper, we prove the following:

\begin{theorem}\label{Theorem 2} Let $\Omega$ be a smooth bounded domain in $\mathbb{R}^4$, $\ell$ be a positive integer, $E_{\ell}^\bot$ and $\lambda_{\ell+1}(\Omega)$ be defined as in (\ref{E}) and (\ref{lambdal}) respectively. Then for any $0\leq\alpha<\lambda_{\ell+1}(\Omega)$, the  supremum
$$\sup_{u\in E_{\ell}^{\bot},\,\| u\|_{2,\alpha}\leq 1}\int_{\Omega}e^{32\pi^2u^2}dx$$
can be attained by some function $u_0\in E_{\ell}^{\bot}\cap C^4(\overline{\Omega}) $ with $\| u_0\|_{2,\alpha}=1$.
\end{theorem}
Obviously Theorem \ref{Theorem 2} extends a result of Yang (\cite{YangJDE-2015}, Theorem 2) and includes (\cite{Nguyen}, Theorem 1.3) as a special case.
The proof of Theorem \ref{Theorem 2} is based on blow-up analysis, which is also used in \cite{Li-Nd,Lu-Yang,Nguyen}.

The remaining part of this paper is organized as follows: In section $2$, we state some results which are crucial in the subsequent analysis; We prove the existence of subcritical maximizers in section $3$ and study the asymptotic behavior of these maximizers in section $4$; In section $5$, we will give upper bound estimates of the functional $\int_{\Omega}e^{32\pi^2u^2}dx$ under the assumption of blow-up analysis; In section $6$, we construct a sequence of test functions to complete the proof of Theorem
{\ref{Theorem 2}}.

\section{Preliminary results}
In this section, we state some preliminaries which would bring a great convenience during our calculation.

Let $G:\Omega\times\Omega\rightarrow \mathbb{R}$ be the Green function of $\Delta^2$ under the Dirichlet condition. That is, for every $y\in\Omega$, the mapping $x\mapsto G(x,y)$ satisfies (in the sense of distribution)
\bna\le\{\begin{array}{lll}\Delta^2 G(x,y)=\delta_{x}(y) \ &{\rm in} \ \Omega,\\[1.2ex]
G(x,y)=\frac{\partial G}{\partial \nu}=0 \ &{\rm on} \ \partial\Omega.
\end{array}\ri.\ena
All functions $u\in W_0^{2,2}(\Omega)\cap C^4(\overline{\Omega})$ satisfying $\Delta^2 u=f(u)$ can be written as
\be\label{u}
u(x)=\int_{\Omega}G(x,y)f(y)dy.
\ee
Now, we collect a property for derivatives of $G$, see for example \cite{D-S}. There exists $C>0$ such that
$$\le|G(x,y)\ri|\leq C\log\le(2+\frac{1}{|x-y|}\ri)$$
and
\be\label{eG2}
|\nabla^kG(x,y)|\leq C|x-y|^{-k},\quad k\geq1 \quad\ee
for all $x,y\in \Omega$, $x\neq y$.
We next recall the Pohozaev identity due to Mitidieri \cite{Miti}.
\begin{lemma}\label{lemma A} Let $\Omega'$ be a smooth bounded domain in $\mathbb{R}^4$, $u\in C^4(\overline{\Omega}')$ be a solution of $\Delta^2u=f(u)$ in $\Omega'$. Then we have for any $y\in \mathbb{R}^4$,
\bna
4\int_{\Omega'}F(u)dx&=&\int_{\partial\Omega'}<x-y,\nu>F(u)d\sigma+\frac{1}{2}\int_{\partial\Omega'}v^2<x-y,\nu>d\sigma+2\int_{\partial\Omega'}\frac{\partial u}{\partial\nu}v d\sigma\\
&&+\int_{\partial\Omega'}\le(\frac{\partial v}{\partial\nu}<x-y,Du>+\frac{\partial u}{\partial\nu}<x-y,Dv>-<Dv,Du><x-y,\nu>\ri)d\sigma
\ena
where $F(u)=\int^u_0 f(t)dt$, $-\Delta u=v$ and $\nu$ is the normal outward derivative of $x$ on $\partial\Omega'$.
\end{lemma}

Similar to \cite{Lu-Yang, Nguyen}, we have the following Lion's type result. Namely
\begin{lemma}\label{lemma B}
Let $\{u_k\}\subset E_{\ell}^{\bot}$ be a sequence of functions and $0\leq\alpha<\lambda_{\ell+1}(\Omega)$ be fixed. If $u_k\rightharpoonup u^*$ weakly in $W_0^{2,2}(\Omega)$ and $\|u_k\|_{2,\alpha}=1$. Then for any $p<(1-\|u^*\|^2_{2,\alpha})^{-1}$
$$\limsup_{k\ra\infty}\int_{\Omega}e^{32\pi^2p u_k^2}dx<\infty.$$
\end{lemma}

\section{Extremals for subcritical Adams inequalities}
In this section, we shall prove that for any $0<\epsilon<32\pi^2$, there exists some function $u_\epsilon\in E_\ell^{\perp}\cap C^4(\overline{\Omega})$ with $\|u_\epsilon\|_{2,\alpha}=1$ such that
\be\label{subc}
\int_{\Omega}e^{(32\pi^2-\epsilon)u_\epsilon^2}dx=\sup_{u\in E_\ell^{\perp},\,\| u\|_{2,\alpha}\leq1}\int_{\Omega}e^{(32\pi^2-\epsilon)u^2}dx
\ee
where $\|\cdot\|_{2,\alpha}$ is defined as in (\ref{1alpha}).

This is based on a direct method in calculus of variations. For any $0<\epsilon<32\pi^2$, we take a sequence of functions $\{u_{k}\}\subset E_\ell^{\perp}$ satisfying that
\be\label{2alpha}\int_{\Omega}|\Delta u_k|^2dx-\alpha\int_{\Omega}u_k^2dx\leq 1\ee
and that
\be\label{ls}\lim_{k\ra\infty}\int_\Omega e^{(32\pi^2-\epsilon)u_k^2}dx=\sup_{u\in E_\ell^{\perp},\,\|u\|_{2,\alpha}\leq 1}\int_{\Omega}e^{(32\pi^2-\epsilon)u^2}dx.\ee
It follows from (\ref{2alpha}) and $0\leq\alpha<\lambda_{\ell+1}(\Omega)$ that $u_{k}$ is bounded in $W_0^{2,2}(\Omega)$. Then there exists some function $u_\epsilon\in W_0^{2,2}(\Omega)$ such that up to a subsequence
\bea\nonumber
&u_k\rightharpoonup u_{\epsilon} &{\rm weakly} \, \ {\rm in}  \, \ W_0^{2,2}(\Omega),\\\nonumber
&u_k\rightarrow u_{\epsilon} &{\rm strongly} \, \ {\rm in} \, \ L^q(\Omega) \, \ (\forall q>1),\\\nonumber
&u_k\ra u_{\epsilon} & a.e.\, \ {\rm in} \, \ \Omega.
\eea
Since $u_k\in E_\ell^{\perp}$, we have
$$ \int_\Omega u_\epsilon e_{ij}dx=\lim_{k\ra \infty}\int_\Omega u_k e_{ij}dx=0,  \quad \forall 1\leq j\leq n_i,\ 1\leq i\leq \ell.\quad$$
Hence $u_\epsilon \in E_\ell^{\perp}$.
By Lemma \ref{lemma B}, we have $e^{(32\pi^2-\epsilon)u_k^2}$ is bounded in $L^r(\Omega)$ for some $r>1$. Therefore
$$\lim_{k\ra\infty}\int_{\Omega}e^{(32\pi^2-\epsilon)u_k^2}dx=\int_{\Omega}e^{(32\pi^2-\epsilon)u_\epsilon^2}dx.$$
This together with (\ref{ls}) immediately leads to (\ref{subc}). Obviously the supremum (\ref{subc}) is strictly greater than the volume of $\Omega$. Thus $u_\epsilon\neq0$. If $\|u_\epsilon\|_{2,\alpha}<1$, we have
\bea\nonumber
\sup_{u\in E_\ell^{\perp},\,\|u\|_{2,\alpha}\leq1}\int_{\Omega}e^{(32\pi^2-\epsilon)u^2}dx&=&\int_{\Omega}e^{(32\pi^2-\epsilon)u_\epsilon^2}dx\\\nonumber
&<&\int_{\Omega}e^{(32\pi^2-\epsilon)u_\epsilon^2/\|u_\epsilon\|^2_{2,\alpha}}dx\\\nonumber
&\leq&\sup_{u\in E_\ell^{\perp},\,\|u\|_{2,\alpha}\leq1}\int_{\Omega}e^{(32\pi^2-\epsilon)u^2}dx,\eea
which is a contradiction.
Then we obtain $\|u_\epsilon\|_{2,\alpha}=1$.

A straightforward calculation shows that $u_\epsilon$ satisfies the following Euler-Lagrange equation:
\be\label{Euler-Lagrange}\le\{\begin{array}{lll}\Delta^2 u_\epsilon-\alpha u_\epsilon=\f{1}{\lambda_\epsilon}u_\epsilon e^{(32\pi^2-\epsilon)u_\epsilon^2}-\sum^{\ell}_{i=1}\sum_{j=1}^{n_i
}\frac{\beta_{ij,\epsilon}}{\lambda_\epsilon}e_{ij} \quad {\rm in} \quad  \Omega,\\[1.2ex]
\int_{\Omega}|\Delta u_\epsilon|^2dx-\alpha\int_{\Omega}u_\epsilon^2dx=1,\\[1.2ex]
u_\epsilon=\frac{\partial u_\epsilon}{\partial \nu}=0 \quad {\rm on} \quad  \partial\Omega,\\[1.2ex]
\lambda_\epsilon=\int_{\Omega} u_\epsilon^2e^{(32\pi^2-\epsilon)u_\epsilon^2}dx,\\[1.2ex]
\beta_{ij,\epsilon}=\int_{\Omega}e_{ij}u_\epsilon e^{(32\pi^2-\epsilon)u_\epsilon^2}dx.
\end{array}\ri.\ee
Applying the standard regularity theory to (\ref{Euler-Lagrange}), we obtain $u_\epsilon\in C^4(\overline{\Omega})$. Since $u_\epsilon$ is bounded in $W_0^{2,2}(\Omega)$, we can assume without loss of generality,
\bea\nonumber
&u_\epsilon\rightharpoonup u_{0} \ &{\rm weakly} \, \ {\rm in} \, \ W_0^{2,2}(\Omega),\\\nonumber
&u_\epsilon\rightarrow u_{0}  \ &{\rm strongly} \, \ {\rm in} \, \ L^q(\Omega) \, \ (\forall q>1),\\\nonumber
&u_\epsilon\ra u_{0}  &a.e. \, \ {\rm in} \, \ \Omega.
\eea
Since $u_\epsilon\in E_\ell^{\perp}$, then
$$ \int_\Omega u_0 e_{ij}dx=\lim_{\epsilon\ra 0}\int_\Omega u_\epsilon e_{ij}dx=0,  \quad 1\leq j\leq n_i, 1\leq i\leq \ell.\quad$$
Hence we get $u_0\in E_\ell^{\perp}$.

Denote $c_\epsilon=|u_\epsilon|(x_\epsilon)=\max_{\Omega}|u_\epsilon|$.
If $c_\epsilon$ is bounded, then $e^{(32\pi^2-\epsilon)u_\epsilon^2}$ is bounded in $L^{\infty}(\Omega)$.
 Clearly $\lambda_\epsilon^{-1}\beta_{ij,\epsilon}e_{ij}$ is also bounded in $L^{\infty}(\Omega)$. Thus for any $u\in E_{\ell}^{\bot}(\Omega)$ with $\| u\|_{2,\alpha}\leq1$, we have by (\ref{subc}) and the Lebesgue dominated convergence theorem
$$\int_{\Omega}e^{32\pi^2u^2}dx=\lim_{\epsilon\ra0}\int_{\Omega}e^{(32\pi^2-\epsilon)u^2}dx\leq\lim_{\epsilon\ra0}\int_{\Omega}e^{(32\pi^2-\epsilon)u_\epsilon^2}dx=\int_{\Omega}e^{32\pi^2u_0^2}dx.$$
This implies that
$$\int_{\Omega}e^{32\pi^2u_0^2}dx=\sup_{u\in E_\ell^{\perp},\,\|u\|_{2,\alpha}\leq 1}\int_{\Omega}e^{32\pi^2u^2}dx.$$
Obviously $\|u_0\|_{2,\alpha}=1$. Using the inequality $e^t\leq1+te^t$ for $t\geq0$, we have
\bea\nonumber\int_{\Omega}e^{(32\pi^2-\epsilon)u_\epsilon^2}dx
&\leq&\int_{\Omega}\le(1+(32\pi^2-\epsilon)u_\epsilon^2e^{(32\pi^2-\epsilon)u_\epsilon^2}\ri)dx\\\nonumber
&=&|\Omega|+(32\pi^2-\epsilon)\lambda_\epsilon.\eea
Hence
$$\liminf_{\epsilon\ra0}\lambda_\epsilon>0.$$Applying the standard regularity theory to (\ref{Euler-Lagrange}), we obtain $u_0\in C^4(\overline{\Omega})$. Therefore $u_0$ is a desired extremal function.

Without loss of generality, we assume there exists some point $x_0\in\overline{\Omega}$ such that
$$x_\epsilon\ra x_0,\quad c_\epsilon=u_\epsilon(x_\epsilon)\ra+\infty \quad as\quad \epsilon\ra0,$$
or we will replace $u_\epsilon$ by $-u_\epsilon$ instead.
In the sequel, we do not distinguish sequence and subsequence, the reader can understand it from the context.

\section{Asymptotic behavior of extremals for subcritical functionals}
In this section, we shall prove that $u_0=0$ and obtain the following Lions type energy concentration result:
\be\label{ec}|\Delta u_\epsilon|^2dx\rightharpoonup\delta_{x_0}\quad {\rm as} \quad\epsilon\rightarrow 0\ee
in the sense of measure, where $\delta_{x_0}$ is the usual Dirac measure centered at $x_0$.

Suppose $u_0\not\equiv0$. In view of Lemma \ref{lemma B}, we have $e^{(32\pi^2-\epsilon)u_\epsilon^2}$ is bounded in $L^r(\Omega)$ for any fixed $r$ with $1<r<1/(1-\|u_0\|_{2,\alpha}^2)$. Note also that $\lambda_\epsilon^{-1}\beta_{ij,\epsilon}e_{ij}$ is bounded in $L^{\infty}(\Omega)$. Applying the standard regularity theory to (\ref{Euler-Lagrange}), we get $u_\epsilon$ is uniformly bounded in $\Omega$, which contradicts with $c_\epsilon\ra\infty$ as $\epsilon\ra0$. Hence $u_0\equiv0$.

Similar to \cite{Lu-Yang, Nguyen}, we can derive $x_0\not\in\partial\Omega$. When $x_0\in \Omega$, suppose (\ref{ec}) is not true. Noting that $\|\Delta u_\epsilon\|_{2,\alpha}=1$, we can find $r_0>0$ and $0<\eta<1$ such that
$$\limsup_{\epsilon\ra0}\int_{\mathbb{B}_{r_0}(x_0)}|\Delta u_\epsilon|^2dx\leq1-\eta.$$
Choose a cut-off function $\phi\in C^2_0(\mathbb{B}_{r_0}(x_0))$, which is equal to $1$ on $\mathbb{B}_{r_0/2}(x_0)$, such that $\mathbb{B}_{r_0}(x_0)\subset \Omega$ and
$$\limsup_{\epsilon\ra0}\int_{\mathbb{B}_{r_0}(x_0)}|\Delta( \phi u_\epsilon)|^2dx\leq1-\eta.$$
By the Adams inequality (\ref{Adams}), $e^{(32\pi^2-\epsilon)\phi^2 u_\epsilon^2}$ is bounded in $L^{r}(\mathbb{B}_{r_0}(x_0))$ for some $r>1$ and thus $e^{(32\pi^2-\epsilon)u_\epsilon^2}$ is bounded in $L^{r}(\mathbb{B}_{r_0/2}(x_0))$ provided that $\epsilon$ is sufficiently small. On the other hand, $\lambda_\epsilon^{-1}\beta_{ij,\epsilon}e_{ij}$ is bounded in $L^{\infty}(\Omega)$. Applying the standard regularity theory to (\ref{Euler-Lagrange}), we derive that $u_\epsilon$ is bounded in $\overline{\mathbb{B}_{r_0/4}(x_0)}$ contradicting $c_\epsilon\ra\infty$. Hence we obtain (\ref{ec}).

Let
$$r_\epsilon^4=\frac{\lambda_\epsilon}{c_\epsilon^2e^{(32\pi^2-\epsilon)c_\epsilon^2}}.$$
Then for any $0<\gamma<32\pi^2$, we have by the H\"older inequality and the Adams inequality (\ref{Adams}),
\be\label{r-epsilon}\lim_{\epsilon\ra0}r_\epsilon^4c_\epsilon^2e^{\gamma
 c_\epsilon^2}=\lim_{\epsilon\ra0}e^{-(32\pi^2-\epsilon-\gamma)c_\epsilon^2}\int_{\Omega}u_\epsilon^2e^{(32\pi^2-\epsilon)u_\epsilon^2}dx=0.\ee
This implies that $r_\epsilon$ converges to zero rapidly. To proceed, we set the following quantities
$$b_\epsilon=\frac{\lambda_\epsilon}{\int_{\Omega}|u_\epsilon|e^{(32\pi^2-\epsilon)u_\epsilon^2}dx},\quad\quad\vartheta=\lim_{\epsilon\rightarrow0}\frac{c_\epsilon}{b_\epsilon},\quad\quad\mu=\lim_{\epsilon\ra0}\frac{\int_{\Omega}u_\epsilon e^{(32\pi^2-\epsilon)u_\epsilon^2}dx}{\int_{\Omega}|u_\epsilon|e^{(32\pi^2-\epsilon)u_\epsilon^2}dx}.$$
Let $\Omega_\epsilon=\{x\in \mathbb{R}^4:x_\epsilon+r_\epsilon x\in \Omega\}.$
Define two blow-up sequences of functions on $\Omega_\epsilon$ as
$$v_\epsilon(x)=\frac{u_{\epsilon}(x_\epsilon+r_\epsilon x)}{c_\epsilon},\quad w_\epsilon(x)=b_\epsilon\le(u_\epsilon(x_\epsilon+r_\epsilon x)-c_\epsilon\ri).$$
A straightforward calculation shows
$$\Delta^2v_\epsilon(x)=\alpha r_\epsilon^4v_\epsilon(x)+\lambda_{\epsilon}^{-1}r_\epsilon^4v_\epsilon(x) e^{(32\pi^2-\epsilon)u_\epsilon^2(x_\epsilon+r_\epsilon x)}-c_\epsilon^{-1}r_\epsilon^4\sum^{\ell}_{i=1}\sum_{j=1}^{n_i
}\frac{\beta_{ij,\epsilon}}{\lambda_\epsilon}e_{ij}\quad {\rm in}\quad\Omega_\epsilon
$$
and
\be\label{ww}
\Delta^2w_\epsilon(x)=\alpha b_\epsilon c_\epsilon r_\epsilon^4v_\epsilon(x)+b_\epsilon c_\epsilon^{-1}v_\epsilon(x)e^{(32\pi^2-\epsilon)c_\epsilon b_\epsilon^{-1}(1+v_\epsilon)w_\epsilon}-b_\epsilon r_\epsilon^4\sum^{\ell}_{i=1}\sum_{j=1}^{n_i
}\frac{\beta_{ij,\epsilon}}{\lambda_\epsilon}e_{ij}\quad {\rm in}\quad\Omega_\epsilon.
\ee
Obviously $|v_\epsilon|\leq1$. Then for any fixed $R>0$ and $x\in\mathbb{B}_{R}(0)$, we obtain
$$\int_{\mathbb{B}_{R}(0)}\le|\Delta v_\epsilon\ri|^2dx=\int_{\mathbb{B}_{Rr_\epsilon}(x_\epsilon)}\frac{\le|\Delta u_\epsilon(y)\ri|^2}{c_\epsilon^2}dy=o_\epsilon(1)$$
and
\be\nonumber
\le|\Delta^2 v_\epsilon(x)\ri|=\le|r_\epsilon^4\le(\alpha v_\epsilon(x)+\lambda_\epsilon^{-1}v_\epsilon e^{(32\pi^2-\epsilon)u_\epsilon^2(x_\epsilon+r_\epsilon x)}-c_\epsilon^{-1}\sum^{\ell}_{i=1}\sum_{j=1}^{n_i
}\frac{\beta_{ij,\epsilon}}{\lambda_\epsilon}e_{ij}\ri)\ri|
=o_\epsilon(1).
\ee
These estimates and the standard regularity theory imply
\be\label{p}
v_{\epsilon}\ra v\quad {\rm in} \quad C^4_{\rm loc}(\mathbb{R}^4).\ee
Since $v(0)=\lim_{\epsilon\ra0}v_\epsilon(0)=1$, we conclude that $v(x)\equiv1$ on $\mathbb{R}^4$ by using the Liouville theorem.

Now we consider the convergence of $w_\epsilon$. Using the Green representation formula (\ref{u}), we get
$$u_\epsilon(x)=\int_{\Omega}G(x,y)\le(\alpha u_\epsilon(y)+\frac{1}{\lambda_\epsilon}u_\epsilon(y) e^{(32\pi^2-\epsilon)u_\epsilon^2(y)}-\sum^{\ell}_{i=1}\sum_{j=1}^{n_i}\frac{\beta_{ij,\epsilon}}{\lambda_\epsilon}e_{ij}\ri)dy.$$
Then for $m=1,\,2$
$$\nabla^{m}u_\epsilon(x)=\int_{\Omega}\nabla_x^mG(x,y)\le(\alpha u_\epsilon(y)+\frac{1}{\lambda_\epsilon}u_\epsilon(y) e^{(32\pi^2-\epsilon)u_\epsilon^2(y)}-\sum^{\ell}_{i=1}\sum_{j=1}^{n_i}\frac{\beta_{ij,\epsilon}}{\lambda_\epsilon}e_{ij}\ri)dy.$$
By (\ref{eG2}) and (\ref{r-epsilon}), we have for any $R>0$ and $x\in \mathbb{B}_R(0)$
\bea\nonumber
\le|\nabla^m w_\epsilon(x)\ri|&=&b_\epsilon r_\epsilon^m\le|\int_{\Omega}\nabla^m_x G(x_\epsilon+r_\epsilon x,y)\Delta^2u_\epsilon(y)dy\ri|\\\nonumber
&\leq&Cb_\epsilon r_\epsilon^m\Bigg(\alpha\int_{\Omega}\frac{|u(y)|}{|x_\epsilon+r_\epsilon x-y|^m}dy+\frac{1}{\lambda_\epsilon}\int_{\mathbb{B}_{2Rr_\epsilon(x_\epsilon)}}\frac{|u_\epsilon(y)|e^{(32\pi^2-\epsilon)u_\epsilon^2(y)}}{|x_\epsilon+r_\epsilon x-y|^m}dy\\\nonumber
&&+\frac{1}{\lambda_\epsilon}\int_{\Omega\backslash\mathbb{B}_{2Rr_\epsilon(x_\epsilon)}}\frac{|u_\epsilon(y)|e^{(32\pi^2-\epsilon)u_\epsilon^2(y)}}{|x_\epsilon+r_\epsilon x-y|^m}dy+\frac{1}{\lambda_\epsilon}\int_{\Omega}\sum^{\ell}_{i=1}\sum_{j=1}^{n_i}\frac{|\beta_{ij,\epsilon}||e_{ij}|}{|x_\epsilon+r_\epsilon x-y|^m}dy\Bigg)\\\nonumber
&\leq&C\Bigg(\alpha b_\epsilon c_\epsilon r_\epsilon^m\int_{\Omega}\frac{1}{|x_\epsilon+r_\epsilon x-y|^m}dy+\frac{b_\epsilon}{c_\epsilon}\int_{\mathbb{B}_{2R}(0)}\frac{1}{|x-z|^m}dz+\frac{1}{R^m}\\\nonumber
&&+\mu r_\epsilon^m\sum^{\ell}_{i=1}\sum_{j=1}^{n_i}\int_{\Omega}\frac{1}{|x_\epsilon+r_\epsilon x-y|^m}dy\Bigg)\\\label{C-1}
&\leq&C(R).
\eea
Here we have used $b_\epsilon\leq c_\epsilon$ and $|\mu|\leq1$. Applying the standard regularity theory to (\ref{ww}), we obtain
\be\label{pp}
w_\epsilon\ra w\quad {\rm in} \quad C^{4}_{\rm loc}(\mathbb{R}^{4}).\ee
If $\vartheta=\lim_{\epsilon\ra0}\frac{c_\epsilon}{b_\epsilon}<+\infty$, then we can see from (\ref{r-epsilon}), (\ref{p}) and (\ref{pp}) that $w$ satisfies
\be\label{r}\le\{\begin{array}{lll}\Delta^2 w(x)=\frac{1}{\vartheta}e^{64\pi^2\vartheta w(x)}\quad{\rm in}\quad \mathbb{R}^4,\\[1.2ex]
w(x)\leq w(0)=0,\\[1.2ex]
\int_{\mathbb{R}^4}e^{64\pi^2\vartheta w(x)}dx<+\infty.
\\[1.2ex]
\end{array}\ri.\ee
To understand $w_\epsilon(x)$ further, we have
$$\Delta w_\epsilon=b_\epsilon r_\epsilon^2\int_{\Omega}\Delta_{x}G(x_\epsilon+r_\epsilon x,y)\le(\alpha u_\epsilon(y)+\frac{1}{\lambda_\epsilon}u_\epsilon(y) e^{(32\pi^2-\epsilon)u_\epsilon^2(y)}-\sum^{\ell}_{i=1}\sum_{j=1}^{n_i}\frac{\beta_{ij,\epsilon}}{\lambda_\epsilon}e_{ij}\ri)dy.$$
Hence for any $R>0$, we obtain by Fubini theorem
\bea\nonumber
\int_{\mathbb{B}_{R}(0)}|\Delta w_\epsilon(x)|dx&\leq&Cb_\epsilon r_\epsilon^2\int_{\Omega}\frac{|u_\epsilon(y)|e^{(32\pi^2-\epsilon)u_\epsilon^2(y)}}{\lambda_\epsilon}\le(\int_{\mathbb{B}_{R}(0)}\frac{1}{|x_\epsilon+r_\epsilon x-y|^2}dx\ri)dy\\\nonumber
&&+C b_\epsilon r_\epsilon^2\int_{\Omega}|u_\epsilon(y)|\le(\int_{\mathbb{B}_{R}(0)}\frac{1}{|x_\epsilon+r_\epsilon x-y|^2}dx\ri)dy\\\nonumber
&&+C b_\epsilon r_\epsilon^2\sum^{\ell}_{i=1}\sum_{j=1}^{n_i}\int_{\Omega}|e_{ij}|\le(\int_{\mathbb{B}_{R}(0)}\frac{1}{|x_\epsilon+r_\epsilon x-y|^2}dx\ri)dy\\\nonumber
&\leq& C R^2.
\eea
This together with (\ref{r}) and the result of \cite{CLin, WeiXu} implies that
$$w(x)=-\frac{1}{16\pi^2\vartheta}\log\le(1+\frac{\pi}{\sqrt{6}}|x|^2\ri),\quad x\in \mathbb{R}^4.$$
When $\vartheta=+\infty$, we have by (\ref{C-1}), $|\Delta w(x)|\leq C R^{-2}$ for $x\in \mathbb{B}_R(0)$.
Letting $R\ra+\infty$, we have $w(x)$ is a harmonic function in $\mathbb{R}^4$. Since $w(x)\leq w(0)=0$, then $w(x)\equiv0$ by the Liouville Theorem.

We next consider the convergence behavior of $u_\epsilon$ away from the blow-up point $x_0$.
Let $\psi_\epsilon$ be a solution of the following Dirichlet problem
\bea\label{G}\le\{\begin{array}{lll}\Delta^2 \psi_\epsilon(x)=\frac{1}{\lambda_\epsilon} b_\epsilon u_\epsilon(x)e^{(32\pi^2-\epsilon)u_\epsilon^2(x)}-\sum^{\ell}_{i=1}\sum_{j=1}^{n_i}\frac{1}{\lambda_\epsilon}b_\epsilon \beta_{ij,\epsilon}e_{ij}\quad{\rm in}\quad \Omega,\\[1.2ex]
\psi_\epsilon(x)=\frac{\partial \psi_\epsilon}{\partial \nu}=0\quad{\rm on}\quad \partial\Omega.
\end{array}\ri.\eea
Using the Green representation formula (\ref{u}), we have
$$\psi_\epsilon(x)=\int_{\Omega}G(x,y)\le(\frac{b_\epsilon}{\lambda_\epsilon} u_\epsilon(y)e^{(32\pi^2-\epsilon)u_\epsilon^2(y)}-\sum^{\ell}_{i=1}\sum_{j=1}^{n_i}\frac{b_\epsilon \beta_{ij,\epsilon}}{\lambda_\epsilon}e_{ij}\ri)dy.$$
By differentiating with respect to $x$ for $m=1,2$, we obtain
$$|\nabla^m \psi_\epsilon(x)|\leq C\le(\int_{\Omega}|x-y|^{-m}\frac{b_\epsilon}{\lambda_\epsilon}|u_\epsilon(y)|e^{(32\pi^2-\epsilon)u_\epsilon^2(y)}dy+\sum^{\ell}_{i=1}\sum_{j=1}^{n_i}\int_{\Omega}|x-y|^{-m}\frac{b_\epsilon}{\lambda_\epsilon}|e_{ij}||\beta_{ij,\epsilon}|dy \ri).$$
For $1<s<2$, applying the basic inequality $(a+b)^s\leq2^{s-1}(a^s+b^s)$ for $a\geq0$ and $b\geq0$ and the H\"older inequality, we obtain
\bea\nonumber
|\nabla^m \psi_\epsilon(x)|^s&\leq&C^s2^{s-1}\le(\int_{\Omega}\frac{b_\epsilon}{\lambda_\epsilon}|u_\epsilon(y)|\frac{e^{(32\pi^2-\epsilon)u_\epsilon^2(y)}}{|x-y|^m}dy\ri)^r+C^s2^{s-1}\le(\sum^{\ell}_{i=1}\sum_{j=1}^{n_i}\int_{\Omega}\frac{b_\epsilon}{\lambda_\epsilon}\frac{|e_{ij}||\beta_{ij,\epsilon}|}{|x-y|^m}dy\ri)^r\\\nonumber
&\leq&C\int_{\Omega}\frac{b_\epsilon}{\lambda_\epsilon}\frac{|u_\epsilon(y)|e^{(32\pi^2-\epsilon)u_\epsilon^2(y)}}{|x-y|^{ms}}dy+C\sum^{\ell}_{i=1}\sum_{j=1}^{n_i}\int_{\Omega}\frac{|e_{ij}|}{|x-y|^{ms}}dy.
\eea
By the Fubini theorem, we have $\|\nabla^m \psi_\epsilon\|_s\leq C$ for $m=1,2$. Hence
\be\label{psi}\|\psi_\epsilon\|_{W_0^{2,s}(\Omega)}\leq C.\ee
Denote $\varphi_\epsilon=b_\epsilon u_\epsilon-\psi_\epsilon$. In view of (\ref{Euler-Lagrange}) and (\ref{G}), we get
\bea\label{var}\le\{\begin{array}{lll}\Delta^2 \varphi_\epsilon(x)=\alpha \le(\varphi_\epsilon(x)+\psi_\epsilon(x)\ri) \ &{\rm in} \,\ \Omega,\\[1.2ex]
\varphi(x)=\frac{\partial \varphi_\epsilon}{\partial \nu}=0 \ &{\rm on} \,\  \partial\Omega.
\end{array}\ri.\eea
Multiplying both sides of (\ref{var}) by $\varphi_\epsilon$, we have by the definition of $\lambda_{\ell+1}(\Omega)$ and the H\"older inequality,
\bna
\int_{\Omega}|\Delta \varphi_\epsilon|^2dx&=&\alpha\int_\Omega \varphi_\epsilon^2dx+\alpha\int_\Omega \psi_\epsilon \varphi_\epsilon dx\\
&\leq&\frac{\alpha}{\lambda_{\ell+1}(\Omega)}\int_{\Omega}|\Delta \varphi_\epsilon|^2dx+\frac{\alpha}{\sqrt{\lambda_{\ell+1}(\Omega)}}\le(\int_{\Omega}\psi_\epsilon^2dx\ri)^{1/2}\le(\int_{\Omega}|\Delta\varphi_\epsilon|^2 dx\ri)^{1/2}.
\ena
Then we get
\be\nonumber\int_{\Omega}|\Delta \varphi_\epsilon|^2dx\leq\frac{\alpha^2 \lambda_{\ell+1}(\Omega)}{(\lambda_{\ell+1}(\Omega)-\alpha)^2}\int_{\Omega}\psi_\epsilon^2 dx.\ee
Hence $\|\varphi_\epsilon\|_{W_0^{2,2}(\Omega)}\leq C$. This together with (\ref{psi}) implies that $b_\epsilon u_\epsilon$ is bounded in $W_0^{2,s}(\Omega)$.
So we can assume there exists some function $G\in W_0^{2,s}(\Omega)$ such that $b_\epsilon u_\epsilon\rightharpoonup G$ weakly in $W_0^{2,s}(\Omega)$ for any $1<s<2$.
Multiplying both sides of (\ref{Euler-Lagrange}) by $b_\epsilon$, we have
\bna\le\{\begin{array}{lll}\Delta^2 (b_\epsilon u_\epsilon)=\alpha (b_\epsilon u_\epsilon)+\f{1}{\lambda_\epsilon}(b_\epsilon u_\epsilon) e^{(32\pi^2-\epsilon)u_\epsilon^2}-\sum^{\ell}_{i=1}\sum_{j=1}^{n_i}\frac{\beta_{ij,\epsilon}}{\lambda_\epsilon}b_\epsilon e_{ij}\quad{\rm in}\quad \Omega,\\[1.2ex]
b_\epsilon u_\epsilon=\frac{\partial (b_\epsilon u_\epsilon)}{\partial \nu}=0\quad{\rm on} \quad \partial\Omega.
\end{array}\ri.\ena
For any fixed $r>0$ such that $\mathbb{B}_{r}(x_0)\subset\Omega$, we derive from the Adams inequality (\ref{Adams}) and the cut-off function theory that $e^{(32\pi^2-\epsilon)u_\epsilon^2}$ is bounded in $L^p(\Omega\backslash\mathbb{B}_{r}(x_0) )$ for some $p>1$. On the other hand, $\lambda^{-1}\beta_{ij,\epsilon}b_\epsilon e_{ij}$ is bounded in $L^{\infty}(\Omega)$. Applying the standard regularity theory we infer that
\be\label{bu}b_\epsilon u_\epsilon\rightarrow G \quad {\rm in} \quad C^4_{\rm loc}(\overline{\Omega}\setminus\{x_0\}).\ee
For any $\phi\in C^\infty(\overline{\Omega})$, we have
\bea\nonumber&&\int_{\Omega}\phi\le(\frac{b_\epsilon u_\epsilon}{\lambda_\epsilon}e^{(32\pi^2-\epsilon)u_{\epsilon}^2}+\alpha b_\epsilon u_\epsilon-\sum^{\ell}_{i=1}\sum_{j=1}^{n_i}\frac{\beta_{ij,\epsilon}}{\lambda_\epsilon}b_\epsilon e_{ij}\ri)dx\\\nonumber
&&=\int_{\Omega}(\phi(x)-\phi(x_0))\frac{b_\epsilon u_\epsilon}{\lambda_\epsilon}e^{(32\pi^2-\epsilon)u_{\epsilon}^2}dx+\int_{\Omega}\phi(x_0)\frac{b_\epsilon u_\epsilon}{\lambda_\epsilon}e^{(32\pi^2-\epsilon)u_{\epsilon}^2}dx\\\nonumber
&&\quad+\alpha\int_{\Omega}\phi b_\epsilon u_\epsilon dx-\sum^{\ell}_{i=1}\sum_{j=1}^{n_i}\int_{\Omega}\frac{\beta_{ij,\epsilon}}{\lambda_\epsilon}\phi b_\epsilon e_{ij} dx
\eea
Lebesgue dominated convergence theorem implies that
\be\label{1}\lim_{\epsilon\ra0}\int_{\{x\in \Omega:,\,|u_\epsilon|\leq1\}}e^{(32\pi^2-\epsilon)u_{\epsilon}^2}dx=|\Omega|.\ee
Combining  (\ref{ls}) and (\ref{1}), we get
\bea\nonumber
\liminf_{\epsilon\ra0}\int_{\Omega}|u_\epsilon|e^{(32\pi^2-\epsilon)u_\epsilon^2}dx&\geq&\liminf_{\epsilon\ra0}\int_{\{x\in\Omega:|u_\epsilon|>1\}}e^{(32\pi^2-\epsilon)u_\epsilon^2}dx\\\nonumber
&=&\sup_{u\in E^+_{\ell},\|u\|_{2,\alpha}\leq1}\int_{\Omega}e^{32\pi^2u^2}dx-|\Omega|\\\nonumber
&>&0.
\eea
This leads to the fact that $b_\epsilon/\lambda_\epsilon$ is bounded. Using the H\"older inequality, we have
\bea\nonumber
&&\int_{\Omega\backslash \mathbb{B}_{r}(x_0)}\le(\phi(x)-\phi(x_0)\ri)\frac{b_\epsilon u_\epsilon}{\lambda_\epsilon}e^{(32\pi^2-\epsilon)u_{\epsilon}^2}dx\\\nonumber
&\leq&C\max_{x\in \overline{\Omega}}\phi(x)\le(\int_{\Omega\backslash \mathbb{B}_{r}(x_0)}e^{(32\pi^2-\epsilon)p_1u_\epsilon^2}dx\ri)^{1/p_1}\le(\int_{\Omega\backslash \mathbb{B}_{r}(x_0)}|u_\epsilon|^{p_2}dx\ri)^{1/p_2}\\\nonumber
&=&o_\epsilon(1)\eea
where $1/p_1+1/p_2=1$. Here we use the facts $u_\epsilon\ra0$ strongly in $L^q(\Omega)$ for any $q>1$ and $e^{(32\pi^2-\epsilon)u_\epsilon^2}$ is bounded in $L^s(\Omega\backslash \mathbb{B}_{r}(x_0))$ for some $s>1$.
By the definition of $b_\epsilon$, we have $$\int_{\mathbb{B}_{r}(x_0)}\frac{1}{\lambda_\epsilon}b_\epsilon|u_\epsilon|e^{(32\pi^2-\epsilon)u_\epsilon^2}dx\leq1.$$
Hence
\be\nonumber
\le|\int_{\mathbb{B}_{r}(x_0)}\le(\phi(x)-\phi(x_0)\ri)\frac{b_\epsilon u_\epsilon}{\lambda_\epsilon}e^{(32\pi^2-\epsilon)u_{\epsilon}^2}dx\ri|\leq
\sup_{x\in \mathbb{B}_{r}(x_0)}\le|\phi(x)-\phi(x_0)\ri|.\ee
We immediately derive
\be\nonumber
\lim_{r\ra 0}\lim_{\epsilon\ra 0}\int_{\Omega}\le(\phi(x)-\phi(x_0)\ri)\frac{b_\epsilon u_\epsilon}{\lambda_\epsilon}e^{(32\pi^2-\epsilon)u_{\epsilon}^2}dx=0.\ee
On the other hand, we can easily get
$$\lim_{\epsilon\ra0}\int_{\Omega}\phi(x_0)\frac{b_\epsilon u_\epsilon}{\lambda_\epsilon}e^{(32\pi^2-\epsilon)u_{\epsilon}^2}dx=\phi(x_0)\lim_{\epsilon\ra0}\frac{\int_{\Omega}u_\epsilon e^{(32\pi^2-\epsilon)u_{\epsilon}^2}dx}{\int_{\Omega}|u_\epsilon | e^{(32\pi^2-\epsilon)u_{\epsilon}^2}dx}=\mu\phi(x_0),$$
$$\lim_{\epsilon\ra0}\int_{\Omega} b_\epsilon u_\epsilon \phi dx=\int_{\Omega} G\phi\ dx,$$
\bna
\lim_{\epsilon\ra0}\le(\sum^{\ell}_{i=1}\sum_{j=1}^{n_i}\int_{\Omega}\frac{\beta_{ij,\epsilon}}{\lambda_\epsilon}b_\epsilon e_{ij}\phi dx\ri)
=\mu\sum^{\ell}_{i=1}\sum_{j=1}^{n_i}e_{ij}(x_0)\int_{\Omega}e_{ij}\phi dx.\ena
These estimates lead to
\bna
&&\int_{\Omega}\phi\le(\frac{b_\epsilon u_\epsilon}{\lambda_\epsilon}e^{(32\pi^2-\epsilon)u_{\epsilon}^2}+\alpha b_\epsilon u_\epsilon-\sum^{\ell}_{i=1}\sum_{j=1}^{n_i}\frac{\beta_{ij,\epsilon}}{\lambda_\epsilon}b_\epsilon e_{ij}\ri)dx\\
&&\quad=\mu\phi(x_0)+\alpha\int_{\Omega}G\phi dx-\mu\sum^{\ell}_{i=1}\sum_{j=1}^{n_i}e_{ij}(x_0)\int_{\Omega}e_{ij}\phi dx.\ena
Therefore, we obtain
\bna\label{G1}\le\{\begin{array}{lll}\Delta^2 G-\alpha G=\mu
\le(\delta_{x_0}-\sum^{\ell}_{i=1}\sum_{j=1}^{n_i}e_{ij}(x_0)e_{ij}\ri)\quad{\rm in}\quad \Omega,\\[1.2ex]
G=\frac{\partial G}{\partial \nu}=0\quad{\rm on} \quad \partial\Omega.
\end{array}\ri.\ena
Take a cut-off function $\eta\in C_0^4(\Omega)$ such that $\eta\equiv1$ on $\mathbb{B}_r(x_0)$ and $\eta\equiv0$ on $\Omega \backslash \mathbb{B}_{2r}(x_0)$, where $\mathbb{B}_{2r}(x_0)\Subset \Omega$. Let
$$g(x)=\frac{\mu}{8\pi^2}\eta(x)\log|x-x_0|+G(x).$$
Then we have
\bna\le\{\begin{array}{lll}\Delta^2 g=f \ &{\rm in}\quad \Omega,\\[1.2ex]
g=\frac{\partial g}{\partial \nu}=0 \ &{\rm on} \quad \partial\Omega
\end{array}\ri.\ena
in a distributional sense, where
\bea\nonumber
f(x)&=&-\frac{\mu}{8\pi^2}\Big(\Delta^2\eta\log|x-x_0|+2\nabla\Delta\eta\nabla\log|x-x_0|+2\Delta\eta\Delta\log|x-x_0|\\\nonumber
&\quad&+2\nabla\eta\nabla\Delta\log|x-x_0|+2\Delta(\nabla\eta\nabla\log|x-x_0|)\Big)+\alpha G-\mu\sum^{\ell}_{i=1}\sum_{j=1}^{n_i}e_{ij}(x_0)e_{ij}.\eea
Note that $f$ is bounded in $L^t(\Omega)$ for any $t>1$. By the standard regularity theory, we obtain $g\in C^3(\overline{\Omega})$. Set $A_{x_0}=g(x_0)$ and
$$\upsilon(x)=g(x)-g(x_0)+\frac{\mu}{8\pi^2}(1-\eta)\log|x-x_0|.$$
Then we get
\be\label{GR}
G(x)=-\frac{\mu}{8\pi^2}\log|x-x_0|+A_{x_0}+\upsilon(x),\ee
where $A_{x_0}$ is a constant depending on $\alpha$ and $x_0$, $\upsilon(x)\in C^3(\overline{\Omega})$ and $\upsilon(x_0)=0$.
Since $u_\epsilon\in E^{\bot}_{\ell}$, then
$$\int_{\Omega}Ge_{ij}dx=\lim_{\epsilon\ra 0}\int_{\Omega}b_\epsilon u_\epsilon e_{ij}dx=0,\quad 1\leq j\leq n_i,\, 1\leq i\leq {\ell}.$$
We have
\be\label{GG}G\in E^{\bot}_{\ell}.\ee

\section{An upper bound}
In this section, we will give an upper bound of the integral $\int_{\Omega}e^{32\pi^2u^2}dx$. The proof is based on the Pohozaev type identity and the capacity estimates.

Set $\Omega'=\mathbb{B}_r(x_\epsilon)$, $y=x_\epsilon$, $u=u_\epsilon$ and $f(u_\epsilon)=\frac{1}{\lambda_\epsilon}u_\epsilon e^{(32\pi^2-\epsilon)u_\epsilon^2}+\alpha u_\epsilon-\sum^{\ell}_{i=1}\sum_{j=1}^{n_i}\frac{ \beta_{ij,\epsilon}}{\lambda_\epsilon}e_{ij}$. Then we have
$$F(u_\epsilon)=\frac{e^{(32\pi^2-\epsilon)u^2_\epsilon}}{2(32\pi^2-\epsilon)\lambda_\epsilon}+\frac{\alpha}{2}u_\epsilon^2-\sum^{\ell}_{i=1}\sum_{j=1}^{n_i}\frac{ \beta_{ij,\epsilon}}{\lambda_\epsilon}e_{ij}u_\epsilon.$$
Applying Lemma \ref{lemma A}, we get for any fixed $r>0$
\bna
\int_{\mathbb{B}_r(x_\epsilon)}e^{(32\pi^2-\epsilon)u^2_\epsilon}dx&=&-\frac{\alpha(32\pi^2-\epsilon)\lambda_\epsilon}{4 b_\epsilon^2}\int_{\mathbb{B}_r(x_\epsilon)}(b_\epsilon u_\epsilon)^2dx+\frac{32\pi^2-\epsilon}{2b_\epsilon}\sum^{\ell}_{i=1}\sum_{j=1}^{n_i}\int_{\mathbb{B}_r(x_\epsilon)}\beta_{ij,\epsilon}e_{ij}b_\epsilon u_\epsilon dx\\
&&+\frac{r}{4}\int_{\partial\mathbb{B}_r(x_\epsilon)}e^{(32\pi^2-\epsilon)u^2_\epsilon}d\sigma+ \frac{\alpha(32\pi^2-\epsilon)\lambda_\epsilon}{4b_\epsilon^2}r\int_{\partial\mathbb{B}_r(x_\epsilon)}(b_\epsilon u_\epsilon)^2d\sigma\\
&&-\frac{32\pi^2-\epsilon}{2b_\epsilon}r\sum^{\ell}_{i=1}\sum_{j=1}^{n_i}\int_{\partial\mathbb{B}_r(x_\epsilon)}\beta_{ij,\epsilon}e_{ij}b_\epsilon u_\epsilon d\sigma\\
&&+\frac{(32\pi^2-\epsilon)\lambda_\epsilon}{4b_\epsilon^2}\le(r\int_{\partial\mathbb{B}_r(x_\epsilon)}|\Delta (b_\epsilon u_\epsilon)|^2d\sigma-4\int_{\partial\mathbb{B}_r(x_\epsilon)}\frac{\partial (b_\epsilon u_\epsilon)}{\partial\nu}\Delta(b_\epsilon u_\epsilon)d\sigma\ri)\\
&&-\frac{(32\pi^2-\epsilon)\lambda_\epsilon}{2b_\epsilon^2}r\int_{\partial\mathbb{B}_r(x_\epsilon)}\le(2\frac{\partial (b_\epsilon u_\epsilon)}{\partial\nu}\frac{\partial (\Delta(b_\epsilon u_\epsilon))}{\partial\nu}-\nabla\Delta(b_\epsilon u_\epsilon)\nabla(b_\epsilon u_\epsilon)\ri)d\sigma
\ena
Letting $\epsilon\ra0$, we have
\bna
\lim_{\epsilon\ra0}\int_{\mathbb{B}_r(x_\epsilon)}e^{(32\pi^2-\epsilon)u^2_\epsilon}dx&=&16\pi^2\lim_{\epsilon\ra0}\frac{\lambda_\epsilon}{b_\epsilon^2}\Bigg(\frac{r}{2}\int_{\partial\mathbb{B}_r(x_\epsilon)}|\Delta G|^2d\sigma-2\int_{\partial\mathbb{B}_r(x_\epsilon)}\frac{\partial G}{\partial\nu}\Delta Gd\sigma\\
&&-2r\int_{\partial\mathbb{B}_r(x_\epsilon)}\frac{\partial G}{\partial \nu}\frac{\partial (\Delta G)}{\partial \nu}d\sigma+r\int_{\partial\mathbb{B}_r(x_\epsilon)}\nabla\Delta G\nabla G d\sigma +o_r(1)+o_{\epsilon,r}(1)\Bigg),\ena
where $o_{\epsilon,r}(1)$ means $\lim_{\epsilon\ra0}o_{\epsilon,r}(1)=0$ for any fixed $r>0$ and $o_{r}(1)$ denotes $\lim_{r\ra0}o_{r}(1)=0$.
By straightforward calculation, we obtain
$$\int_{\partial\mathbb{B}_r(x_0)}|\Delta G|^2d\sigma=\frac{\mu^2}{8\pi^2r}+o(1),\quad\quad\int_{\partial\mathbb{B}_r(x_0)}\frac{\partial G}{\partial \nu}\Delta Gd\sigma=\frac{\mu^2}{16\pi^2}+o(1),$$
$$\int_{\partial\mathbb{B}_r(x_0)}\frac{\partial G}{\partial \nu}\frac{\Delta G}{\partial \nu}d\sigma=-\frac{\mu^2}{8\pi^2r}+o(1),\quad\quad \int_{\partial\mathbb{B}_r(x_0)}\nabla\Delta G \nabla Gd\sigma=-\frac{\mu^2}{8\pi^2r}+o(1).$$
Therefore, we have
$$\lim_{r\ra 0}\lim_{\epsilon\ra 0}\int_{\mathbb{B}_r(x_0)}e^{(32\pi^2-\epsilon)u_\epsilon^2}dx=\mu^2\lim_{\epsilon\ra 0}\frac{\lambda_\epsilon}{b_\epsilon^2}.$$
Since $|\Delta u_\epsilon|^2\rightharpoonup \delta_{x_0}$ in the sense of measure, we get
$$\lim_{r\ra 0}\lim_{\epsilon\ra 0}\int_{\Omega\backslash\mathbb{B}_r(x_0)}e^{(32\pi^2-\epsilon)u_\epsilon^2}dx=|\Omega|.$$
By these two identities we have
$$\lim_{\epsilon\ra 0}\int_{\Omega}e^{(32\pi^2-\epsilon)u_\epsilon^2}dx=|\Omega|+\mu^2\lim_{\epsilon\ra 0}\frac{\lambda_\epsilon}{b_\epsilon^2}$$
provided that $\mu^2>0$.
Similar as that in (\cite{Lu-Yang}, Lemma 4.6), we have $\mu=1$. The proof is omitted. Hence, (\ref{GR})
can be restated as
\bna\label{G1}\le\{\begin{array}{lll}\Delta^2 G-\alpha G=
\delta_{x_0}-\sum^{\ell}_{i=1}\sum_{j=1}^{n_i}e_{ij}(x_0)e_{ij}\quad{\rm in}\quad \Omega,\\[1.2ex]
G=\frac{\partial G}{\partial \nu}=0\quad{\rm on} \quad \partial\Omega.
\end{array}\ri.\ena
Moreover, $G$ can be represented as
\bna\label{G2}
G(x)=-\frac{1}{8\pi^2}\log|x-x_0|+A_{x_0}+\upsilon(x),\ena
where $A_{x_0}$ is a constant depending on $x_0$ and $\alpha$, $\upsilon(x)\in C^3(\overline{\Omega})$ and $\upsilon(x_0)=0$.

We now point out the following results:
$$\lim_{\epsilon\ra0}\vartheta=\lim_{\epsilon\ra 0}\frac{c_\epsilon}{b_\epsilon}=1$$
and $\lambda_\epsilon/c_\epsilon^2$ is bounded.
The reader can refer to \cite{Lu-Yang, Nguyen}  for details.
In particular, $b_\epsilon$ can be replaced by $c_\epsilon$ in (\ref{bu}). Namely, $c_\epsilon u_\epsilon\rightarrow G$ in $C^4_{\rm loc}(\Omega\backslash\{x_0\})$.

The technique of capacity estimates was first applied to deal with first derivatives of Moser inequality \cite{Li2001}. Slightly modified the proof in (\cite{Nguyen}, Section $4$) which is adapted from the idea of (\cite{Lu-Yang}, Section $5$), we have
\begin{lemma}\label{Lemma 2} For any $u\in E^\bot_\ell$ with $\|u\|_{2,\alpha}\leq1$, there holds
\be\label{ub}\sup_{u\in E^\bot_\ell,\,\|u\|_{2,\alpha}\leq1}\int_{\Omega}e^{32\pi^2u^2}dx\leq|\Omega|+\frac{\pi^2}{6}e^{\frac{5}{3}+32\pi^2A_{x_0}}.\ee
\end{lemma}

\section{Test function computation}
In this section, we will construct a sequence of text functions $\phi^*_\epsilon\in E^{\bot}_{\ell}$ such that ${\|\phi^*_\epsilon\|}_{2,\alpha}=1$ and
\be\label{C}
\int_{\Omega}e^{32\pi^2{\phi^*_\epsilon}^{2}}dx>|\Omega|+\frac{\pi^2}{6}e^{\frac{5}{3}+32\pi^2A_{x_0}}
\ee
for $\epsilon>0$ sufficiently small. This leads to a contraction with (\ref{ub}). Hence, blow-up can not occur and thus $c_\epsilon$ must be bounded. The standard regularity theory leads to the existence of the desired extremal function. The proof of Theorem \ref{Theorem 2} is completely finished.

To prove (\ref{C}), we write $r=|x-x_0|$. Recall that $G(x)=-\frac{1}{8\pi^2}\log|x-x_0|+A_{x_0}+\upsilon(x)$. Set
$$\phi_\epsilon=\le\{
     \begin{array}{llll}
     &c+\f{a-\frac{1}{16\pi^2}\log\le(1+\frac{\pi r^2}{\sqrt{6}\epsilon^2}\ri)+A_{x_0}+\upsilon+br^2}{c},
     \quad &r\leq R\epsilon,\\[1.2ex]
     &\f{G}{c},\quad & r>R\epsilon,
     \end{array}
     \ri.$$
where $a$, $b$ and $c$ are constants of $\epsilon$ to be determined later, $R=-\log\epsilon$. In order to assure that $\phi_\epsilon\in W^{2,2}_{0}(\Omega)$, we require
$$\lim_{r\rightarrow(R\epsilon)^-}\phi_\epsilon=\lim_{r\rightarrow(R\epsilon)^+}\phi_\epsilon$$
and
$$\lim_{r\rightarrow(R\epsilon)^-}\nabla\phi_\epsilon=\lim_{r\rightarrow(R\epsilon)^+}\nabla\phi_\epsilon.$$
Then we have
\bna\label{ab}\le\{\begin{array}{lll}a=-c^2-\frac{\log(R\epsilon)}{8\pi^2}+\frac{\log\le(1+\frac{\pi}{\sqrt{6}R^2}\ri)}{16\pi^2}-\frac{1}{16\pi^2\le(1+\frac{\pi}{\sqrt{6}}R^2\ri)} ,\\[1.2ex]
b=-\frac{1}{16\pi^2R^2\epsilon^2\le(1+\frac{\pi}{\sqrt{6}}R^2\ri)}.\\[1.2ex]
\end{array}\ri.\ena
A straightforward calculation shows
$$\int_{\Omega}|\Delta\phi_\epsilon|^2dx-\alpha\int_{\Omega}\phi_\epsilon^2dx=\frac{1}{c^2}\le(-\frac{1}{8\pi^2}\log\epsilon+\frac{1}{16\pi^2}\log\frac{\pi}{\sqrt{6}}+A_{x_0}-\frac{5}{96\pi^2}+O\le(\frac{1}{\log^2 \epsilon}\ri)\ri).$$
Setting $\|\phi_\epsilon\|_{2,\alpha}=1$, we obtain
$$c^2=-\frac{\log\epsilon}{8\pi^2}+\frac{\log\pi}{16\pi^2}-\frac{\log6}{32\pi^2}-\frac{5}{96\pi^2}+A_{x_{0}}+O\le(\frac{1}{\log^2\epsilon}\ri).$$
We calculate on $\mathbb{B}_{R\epsilon}(x_0)$
$$32\pi^2\phi_\epsilon^2\geq\log\frac{\pi^2}{6\epsilon^4}-4\log\le(1+\frac{\pi r^2}{\sqrt{6}\epsilon^2}\ri)+\frac{5}{3}+32\pi^2A_{x_0}+O\le(\frac{1}{\log^2\epsilon}\ri).$$
Hence
\bea\nonumber
\int_{\mathbb{B}_{R\epsilon}(x_0)}e^{32\pi^2\phi_\epsilon^2}dx&=&\frac{\pi^2}{6\epsilon^4}e^{\frac{5}{3}+32\pi^2A_{x_0}}\int_{\mathbb{B}_{R\epsilon}(x_0)}\le(1+\frac{\pi}{\sqrt{6}}\frac{r^2}{\epsilon^2}\ri)^{-4}dx\\\label{B1}
&=&\frac{\pi^2}{6}e^{\frac{5}{3}+32\pi^2A_{x_0}}\le(1+O\le(\frac{1}{\log^2\epsilon}\ri)\ri).\eea
On the other hand, we get on $\Omega\setminus\mathbb{B}_{R\epsilon}(x_0)$
\bea\nonumber
\int_{\Omega\setminus\mathbb{B}_{R\epsilon}(x_0)}e^{32\pi^2\phi_\epsilon^2}dx&\geq&\int_{\Omega\setminus\mathbb{B}_{R\epsilon}(x_0)}(1+32
\pi^2\phi^2_\epsilon)dx\\\label{B2}
&=&|\Omega|+32\pi^2\frac{\|G\|^2_2}{c^2}+O\le(\frac{1}{\log^2\epsilon}\ri).\eea
Combining (\ref{B1}) and (\ref{B2}), we conclude
\be\label{ed}\int_{\Omega}e^{32\pi^2\phi_\epsilon^2}dx\geq|\Omega|+\frac{\pi^2}{6}e^{\frac{5}{3}+32\pi^2A_{x_0}}+\frac{32\pi^2}{c^2}\|G\|^2_2+O\le(\frac{1}{\log^2\epsilon}\ri).\ee
Recalling that $(e_{ij})$ $(1\leq j\leq n_i,\,1\leq i\leq \ell)$ is a basis  of $E_{\ell}^{\bot}$ verifying (\ref{ei}), we set
$$\widetilde{\phi}_\epsilon=\phi_\epsilon-\sum^{\ell}_{i=1}\sum_{j=1}^{n_i}(\phi_\epsilon, e_{ij})e_{ij},$$
where
$$(\phi_\epsilon, e_{ij})=\int_{\Omega}\phi_\epsilon e_{ij}dx.$$
Obviously $\widetilde{\phi}_\epsilon\in E_{\ell}^{\bot}$. A straightforward calculation gives
\bea\nonumber
(\phi_\epsilon, e_{ij})&=&\int_{\mathbb{B}_{R_\epsilon}(x_0)}\le(c+\f{a-\frac{1}{16\pi^2}\log\le(1+\frac{\pi r^2}{\sqrt{6}\epsilon^2}\ri)+A_{x_0}+\upsilon+br^2}{c}\ri)e_{ij}dx\\\label{eij}
&&+\int_{\Omega \backslash \mathbb{B}_{R_\epsilon}(x_0)}\frac{G}{c}e_{ij}dx=o\le(\frac{1}{\log^2\epsilon}\ri).\eea
Here we have used (\ref{GG}) to obtain
$$\int_{\Omega\backslash \mathbb{B}_{R_\epsilon}(x_0)}\frac{G}{c}e_{ij}dx=-\int_{\mathbb{B}_{R_\epsilon}(x_0)}\frac{G}{c}e_{ij}dx=O(\epsilon^4(-\log \epsilon)^{9/2})=o\le(\frac{1}{\log ^2\epsilon}\ri).$$
In view of (\ref{eij}) and the fact $\|\phi_\epsilon\|^2_{2,\alpha}=1$, we obtain
\be\label{phi1}
\widetilde{\phi}_\epsilon=\phi_\epsilon+o\le(\frac{1}{\log^2\epsilon}\ri),\ee
\be\label{phi2}
\|\widetilde{\phi}_\epsilon\|_{2,\alpha}^2=1+o\le(\frac{1}{\log^2\epsilon}\ri).\ee
Combining (\ref{ed}), (\ref{phi1}) and (\ref{phi2}), we derive
\bna
\int_{\Omega}e^{32\pi^2\frac{\widetilde{\phi}_\epsilon^2}{\|\widetilde{\phi}_\epsilon\|_{2,\alpha}^2}}dx&=&\int_{\Omega}e^{32\pi^2\phi^2_\epsilon+o(\frac{1}{\log\epsilon})}dx\\
&\geq&\le(1+o\le(\frac{1}{\log\epsilon}\ri)\ri)\le(|\Omega|+\frac{\pi^2}{6}e^{\frac{5}{3}+32\pi^2A_{x_0}}+\frac{32\pi^2}{c^2}\le(\|G\|^2_2+o(1)\ri)\ri)\\
&\geq&|\Omega|+\frac{\pi^2}{6}e^{\frac{5}{3}+32\pi^2A_{x_0}}+32\pi^2\frac{\|G\|^2_2}{c^2}+o\le(\frac{1}{c^2}\ri).
\ena
Set $\phi_\epsilon^*={\widetilde{\phi}_\epsilon}/{\|\widetilde{\phi}_\epsilon\|_{2,\alpha}}.$
 Noting that $\widetilde{\phi}_\epsilon\in E_\ell^{\bot}$, we get $\phi_\epsilon^*\in E_\ell^{\bot}$. Moreover $\|\phi_\epsilon^*\|_{2,\alpha}=1$ and (\ref{C}) holds.
 The contradiction between (\ref{ub}) and (\ref{C}) shows that $c_\epsilon$ must be bounded. This completes the proof of Theorem {\ref{Theorem 2}}.

\bigskip

{\bf Acknowledgements}. This work was partly supported by the Natural Science Foundation of the Education Department of Anhui Province \,(KJ2016A641) and the Outstanding Young Talents Program of the Education Department of Anhui Province (gxyq2018160).

\end{document}